\documentclass[11pt]{article}
\usepackage{graphicx}
\usepackage{amssymb}

\newtheorem{theorem}{Theorem}[section]
\newtheorem{example}[theorem]{Example}
\newtheorem{lemma}[theorem]{Lemma}

\newtheorem{proposition}[theorem]{Proposition}
\newtheorem{corollary}[theorem]{Corollary}
\title{\textbf{MASS IN THE HYPERBOLIC PLANE}\footnote{This research was supported in part by University of Kansas General Research Allocation 2301559-003}}
\author{Saul Stahl}

\begin{document}
\maketitle
\begin{flushright}
Department
of Mathematics\\
University of Kansas\\
Lawrence,
KS 66045, USA\\
\textbf{stahl@math.ku.edu}
\end{flushright}

\section{INTRODUCTION}

Archimedes computed the center of mass of several regions and bodies [Dijksterhuis], and this fundamental physical notion may very well be due to him.  He based his investigations of this concept on the notion of moment as it is used in his Law of the Lever.  A hyperbolic version of this law was formulated in the nineteenth century leading to the notion of a hyperbolic center of mass of two point-masses [Andrade, Bonola]. In 1987 Galperin proposed an axiomatic definition of the center of mass of finite systems of point-masses in Euclidean, hyperbolic and elliptic $n$-dimensional spaces and proved its uniqueness. His proof is based on Minkowskian, or relativistic, models and evades the issue of moment.  A surprising aspect of this work is that hyperbolic mass is not additive.   Ungar [2004] used the theory of gyrogroups to show that in hyperbolic geometry the center of mass of three point-masses of equal mass coincides with the point of intersection of the medians.  Some information regarding the centroids of finite point sets in spherical spaces can be found in [Fog, Fabricius-Bjerre].

In this article we offer a physical motivation for the hyperbolic Law of the Lever and go on to provide a model-free definition and development of the notions of center of mass, moment, balance and mass of finite point-mass systems in hyperbolic geometry. All these notions are then extended to linear sets and laminae.  Not surprisingly, the center of mass of the uniformly dense hyperbolic triangle coincides with the intersection of the triangle's medians.  However, it is pleasing that a hyperbolic analog of Archimedes's mechanical method can be brought to bear on this problem. The masses of uniform disks and regular polygons are computed in the Gauss model and  these formulas are very surprising.  Other configurations are examined as well.

For general information regarding the hyperbolic plane the reader is referred to [Greenberg, Stahl]

\section{THE HYPERBOLIC LAW OF THE LEVER}

Many hyperbolic formulas can be obtained from their Euclidean analogs by the mere replacement of a length $d$ by $\sinh d$. The Law of Sines and the Theorems of Menelaus and Ceva (see Appendix) are cases in point.  It therefore would make sense that for a lever in the hyperbolic plane a suitable definition of the \emph{moment} of a force $w$ acting perpendicularly at distance $d$ from the fulcrum is
\[w\sinh d\]
Nevertheless, a more physical motivation is in order. We begin with an examination of the balanced weightless lever of Figure 1. This lever is pivoted at $E$ and has masses of weights $w_1$ and $w_2$ at $A$ and $B$ respectively.  By this is meant that there is a mass $D$, off the lever, which exerts attractive forces $\vec{w}_1$ and $\vec{w}_2$ along the straight lines $AD$ and $BD$.  Since this system is assumed to be in equilibrium, it follows that the resultant of the forces $\vec{w}_1$ and $\vec{w}_2$ acts along the straight line $ED$.  Neither the direction  nor the intensity of the resultant are affected by the addition of a pair of equal but opposite forces $\vec{f}_1$ and $\vec{f}_2$ at $A$ and $B$. (Here and below we employ the convention that the magnitude of the vector $\vec{v}$ is denoted by $v$.) We assume that the common magnitude of $f_1$ and $f_2$ is large enough so that the lines of direction  of the partial resultants $\vec{r}_i = \vec{f}_i + \vec{w}_i$, i = 1,2,  intersect in some point, say $C$.  Note that the quadrilateral $ACBD$ lies in the hyperbolic plane whereas the parallelograms of forces at $A$ and $B$ lie in the respective Euclidean tangent planes. This is the standard operating procedure in mathematical physics.

\begin{figure}
\center
\includegraphics{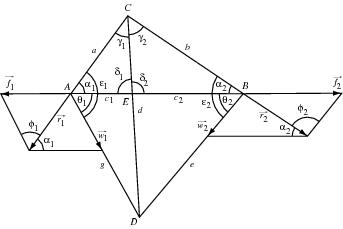}
\caption{}
\end{figure}

 It is now demonstrated that such a system in equilibrium must satisfy the equation
\begin{equation}F_1\sinh c_1 = F_2\sinh c_2
\end{equation}
where each $\vec{F}_i$ is the component of $\vec{w}_i$ in the direction orthogonal to $AB$.
Indeed, it follows from several applications of both the Euclidean and the hyperbolic Laws of Sines that
\[\frac{w_1\sinh c_1}{w_2\sinh c_2} = 
\frac{w_1 \sin \gamma_1 \cdot \frac{\sinh a}{\sin \delta_1}}
{w_2 \sin \gamma_2 \cdot \frac{\sinh b}{\sin \delta_2}}=
\frac{w_1 \sin \gamma_1 \sinh a}{w_2 \sin \gamma_2 \sinh b} =
\]
\[\frac{\sin \gamma_1 \cdot \frac{w_1}{\sin \alpha_1}}
	{\sin \gamma_2 \cdot \frac{w_2}{\sin \alpha_2}} =
	 \frac{\sin\gamma_1 \cdot \frac{ f_1}{\sin \phi_1}}{\sin\gamma_1 \cdot \frac{ f_2}{\sin\phi_2}} = \frac{\sin \gamma_1}{\sin \gamma_2}\  \frac{\sin\epsilon_2}{\sin \epsilon_1} =
\]
\[\frac{\frac{\sin \epsilon_2}{\sin \gamma_2}}{\frac{\sin \epsilon_1}{\sin \gamma_1}} = \frac{\frac{\sinh d }{\sinh e}}{\frac{\sinh d}{\sinh g}} =
\frac{\sinh g}{\sinh e} = 
\frac{\sin\theta_2}{\sin\theta_1}
\]
and Eq'n (1) follows by cross-multiplication.


If we take the mass at $D$ out of the picture and stipulate that $\vec{F_1}$ and $\vec{F_2}$ are simply two forces that act perpendicularly to the lever $AB$ (Fig. 2) then it is makes sense to regard the quantities 
\[F_1\sinh c_1 \ \ \emph{\emph{and}} \ \ F_2\sinh c_2
\]
as the respective moments of the forces $\vec{F_1}$ and $\vec{F_2}$ with respect to the pivot point $E$.  This facilitates the derivation of the resultant of $\vec{F_1}$ and $\vec{F_2}$.
\begin{figure}
\center
\includegraphics{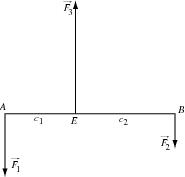}
\caption{}
\end{figure}
 Suppose $c_1$, $c_2$ and $\vec{F}_3 \perp AB$ are  such that
\begin{equation}F_1\sinh c_1  = F_2\sinh c_2  \ \ \ \emph{\emph{and}}\ \ \ \ F_3 = F_1\cosh c_1  + F_2\cosh c_2 
\end{equation}
Then the moments of $\vec{F}_3$ with respect to $A$ and $B$ are, respectively
\[(F_1\cosh c_1  + F_2\cosh c_2 )\sinh c_1  \] 

\[ = F_2\cosh c_1 \sinh c_2  + F_2\sinh c_1 \cosh c_2  = F_2\sinh(c_1 + c_2)
\]
and 
\[(F_1\cosh c_1  + F_2\cosh c_2 )\sinh c_2  = F_1\sinh(c_1 + c_2).
\]Since the right hand sides of these two equations, are, respectively, the moments of $\vec{F}_2$ with respect to $A$ and the moment of $\vec{F}_1$ with respect to $B$, it follows that the equations of (2) do indeed imply equilibrium.  Consequently, the reverse of $\vec{F}_3$ is indeed the resultant of $\vec{F}_1$ and $\vec{F}_2$.
\section{
FINITE POINT-MASS SYSTEMS}

The physical considerations of the previous section motivate the following formal definitions.  A \emph{point-mass} is an ordered pair $(X, x)$ where its \emph{location} $X$ is a point of the hyperbolic plane and its \emph{weight} $x$ is a positive real number.  The 
 \emph{(unsigned) moment} of the point-mass $(X, x)$ with respect to the point $N$ or the straight line $n$ is, respectively, 
\[  M_N(X, x) = x \sinh d(X, N) \ \ \ \  \emph{\emph{or}}\ \ \ \  M_n(X, x) = x \sinh d(X,n)  
\] where  $d(X, N)$ and $d(X, n)$ are the respective hyperbolic distances from  $X$ to $N$ and $n$.

Given any two point-masses $(X,x)$ and $(Y,y)$, their \emph{center of mass} or \emph{centroid} $(X, x)*(Y, y)$ is the point-mass $(Z, z)$, where $Z$ is that point between $X$ and $Y$ such that
\[x\sinh XZ  = y \sinh YZ
\]and

\begin{equation}z = x\cosh XZ  + y \cosh YZ 
\end{equation}Note that this means that the two point-masses have equal moments with respect to their centroid. Moreover, if $X = Y$ then $(X, x)*(Y, y) = (X, x+y)$.

The next two propositions demonstrate that the center of mass "balances" its two constituent point-masses.
\begin{proposition}If $(Z,z) = (X, x) * (Y, y)$ then any two of these point-masses have equal moments with respect to the location of the third one.
\end{proposition}
Proof: It follows from the definitions and that $(X,x)$ and $(Y,y)$ have equal moments with respect to $Z$.  Hence it only remains to show that $(X, x)$ and $(Z, z)$ have equal moments with respect to $Y$.  In other words, that
\[x\sinh(XZ + YZ) = z\sinh YZ 
\]or
\[x\sinh XZ \cosh YZ  + x\cosh XZ \sinh YZ  \]\[ = x\cosh XZ \sinh YZ +y\cosh YZ \sinh YZ 
\]and this equation follows from the fact that $(Z, z)$ is the centroid of $(X, x)$ and  $(Y, y)$.   \hfill Q.E.D.
\newline
Given any two point-masses $(X, x)$ and $(Y, y)$, their \emph{external centroid} is the point-mass $(Z, z)$ such that $Z$ is on the straight line $XY$ but outside the segment joining $X$ and $Y$,
\[x\sinh XZ  = y \sinh YZ
\]and

\[z = |x\cosh XZ - y\cosh YZ|
\]
\begin{figure}
\center
\includegraphics{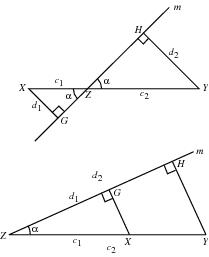}
\caption{}
\end{figure}
\begin{proposition} The two point-masses $(X, x)$ and $(Y, y)$ have equal moments with respect to the intersecting straight line $m$, if and only if $m$ contains at least one of their centroids.  They have equal moments with respect to \underline{every} straight line if and only if they are identical.
\end{proposition}
PROOF: In both of the the diagrams of Figure 3
\[\frac{\sinh d_1}{\sinh c_1} = \frac{\sinh d_2}{\sinh c_2} (= \sin \alpha)
\]and consequently
\[\frac{\sinh d_1}{\sinh d_2} = \frac{\sinh c_1}{\sinh c_2}.
\]
Hence
\[x\sinh d_1 = y \sinh d_2 \ \ \ \ \emph{\emph{if and only if}}\ \ \ \ x\sinh c_1 = y\sinh c_2.\] This implies the first half of the proposition.  The second half follows immediately from the first one.  \hfill Q.E.D.\newline

The following proposition implies that the \emph{center of mass} $\mathcal{C(X)} = (C, c)$ of any finite point-mass system $\mathcal{X}$ is well defined.  This definition clearly satisfies the axioms of [Galperin] and so the two are equivalent.  

\begin{proposition}The binary operation $"*"$ is both commutative and associative.
\end{proposition}
Proof: The commutativity of $''*''$  follows immediately from its definition.  To prove its associativity, let $(X, x), (Y, y), (Z, z)$ be three arbitrary point-masses, and let $(P, p) = (Y, y)*(Z, z), (Q, q) = (Z, z)*(X, x), (R , r) = (X, x)*(Y, y)$ (Fig. 4).   We may assume that $X, Y,$ and $Z$ are not collinear since the degenerate cases follow by an easy independent argument or can be verified from the assumed case by a continuity argument.  Then
\[1 = \frac{x}{y}\cdot\frac{y}{z} \cdot \frac{z}{x} = \frac{\sinh b_1 }{\sinh a_2 }\cdot \frac{\sinh c_1 }{\sinh b_2 }\cdot \frac{\sinh a_1 }{\sinh c_2 }
\]
and the hyperbolic Theorem of Ceva implies that the cevians $XP, YQ, ZR$ are concurrent, say at $T$.

\begin{figure}
\center
\includegraphics{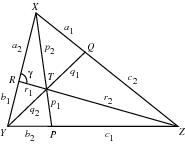}
\caption{}
\end{figure}

Next we show that the point-masses $(R, r)$ and $(Z,z)$ have equal moments with respect to $T$.   In other words, that
\[\frac{\sinh r_1 }{\sinh r_2 } = \frac{z}{x\cosh a_2  +y \cosh b_1 }
\]
However, an application of the unsigned version of the hyperbolic Theorem of Menelaus to $\Delta RYZ$ yields,
\[\frac{\sinh r_1 }{\sinh r_2 } = \frac{\sinh b_2 }{\sinh c_1 }\cdot\frac{\sinh a_2 }{\sinh(a_2 + b_1)}
\]
and hence it suffices to prove that
\[x\cosh a_2 \sinh b_2 \sinh a_2  + y\cosh b_1 \sinh b_2 \sinh a_2  
\]
\[ = z\sinh c_1 \sinh(a_2 + b_1)
\]
This, however, follows easily from the substitutions 
\[x\sinh a_2  = y\sinh b_1 
\]
\[z\sinh c_1  = y\sinh b_2 
\]
and the formula for $\sinh(\alpha + \beta)$.  

This shows that
\begin{equation}[(X, x) * (Y, y)]*(Z, z)\end{equation} is located at $T$.  Because of the symmetry of the construction of $T$ it may be concluded that the same holds for every one of the systems obtained by permuting the constituents of Eq'n (4).

Finally, note that if
\[[(X, x)*(Y, y)]*(Z, z) = (T, t),\]
then,  by several   applications of the hyperbolic Law of Cosines and Eq'n (3),
\[t = x\cosh a_2 \cosh r_1  + y\cosh b_1 \cosh r_1  + z\cosh r_2 
\] 
\[= x(\cosh p_2 - \cos\gamma \sinh a_2 \sinh r_1)\]
\[ + y[\cosh q_2 - \cos (\pi - \gamma) \sinh b_1 \sinh r_1] + z\cosh r_2)
\]
\[= x\cosh p_2 + y\cosh q_2 + z\cosh r_2\]
The pleasing symmetry of this expression demonstrates that all the permutations of  (4)  also have the same masses.
\hfill Q.E.D. \newline

In contrast with masses, moments are additive in the following sense.

\begin{theorem}
Let (X, x) and (Y, y) be two point-masses, and let $m$ be any directed straight line then
\[M_m\left((X, x)*(Y, y)\right) = M_m(X, x) + M_m(Y, y)
\]
\end{theorem}
Proof:  As this is trivial when $X = Y$ we assume that $X$ and $Y$ are distinct.  We first suppose that $XY || m$.  In that case they are known to have a common perpendicular line, say $p$ (Fig. 5).
\begin{figure}
\center
\includegraphics{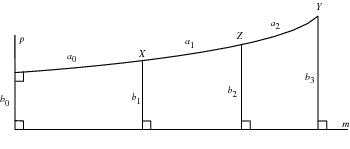}
\caption{}
\end{figure}
By Equation (i) on p. 344 of [Greenberg]
\[\sinh b_1  = \cosh a_0 \sinh b_0
\]
\[\sinh b_2  = \cosh(a_0+a_1)\sinh b_0
\]
\[\sinh b_3 = \cosh(a_0+a_1+a_2)\sinh b_0
\]The proposed equation is now proved by observing that each of the following equations is equivalent to the next.
\[\left(x\cosh a_1  + y\cosh a_2 \right)\sinh b_2  = x\sinh b_1  + y \sinh b_3
\]

\[\left(x\cosh a_1  + y\cosh a_2 \right)\cosh(a_0 + a_1) = x\cosh a_0  + y\cosh (a_0+ a_1 + a_2) \]

\[x\cosh^2 a_1 \cosh a_0  + x\cosh a_1 \sinh a_0 \sinh a_1 
\]
\[+ y\cosh a_2 \cosh a_0 \cosh a_1  + y\cosh a_2 \sinh a_0 \sinh a_1 
\]
\[= x\cosh a_0  + y\cosh a_0 \cosh a_1 \cosh a_2  + y\sinh a_0 \sinh a_1 \cosh a_2  
\]
\[+ y\sinh a_0 \cosh a_1 \sinh a_2  + 
y\cosh a_0 \sinh a_1 \sinh a_2 
\]

\[x\cosh a_0 + x\cosh a_0\sinh^2 a_1 + x\sinh a_0 \sinh a_1 \cosh a_1 
\]
\[= x\cosh a_0 + x\sinh a_0 \sinh a_1 \cosh a_1  +   x\cosh a_0\sinh^2a_1 
\]

\[0 = 0.
\]
\begin{figure}
\center
\includegraphics{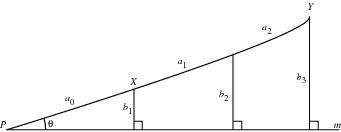}
\caption{}
\end{figure}
On the other hand, if $XY$ and $m$ intersect, say at $P$ (Fig. 6) with $X$ and $Y$ on the same side of $m$, then, by Theorem 8.4ii of [Stahl], 

\begin{equation}\sin\theta = \frac{\sinh b_1 }{\sinh a_0} = \frac{\sinh b_2 }{\sinh(a_0 + a_1)} = \frac{\sinh b_3}{\sinh(a_0 + a_1 + a_2)}
\end{equation}

The required equation is tantamount to 
\[\left(x\cosh a_1  + y\cosh a_2 \right)\sinh b_2  = x\sinh b_1  +  y\sinh b_3
\]which by Eq'n (5) is tantamount to 

\[\left(x\cosh a_1  + y\cosh a_2 \right)\sinh(a_0+a_1) \]\[= x\sinh a_0 +  y\sinh(a_0+a_1+a_2)\]
This however, is easily proved by the same technique as was used in the first half of this proof.

If  $X$ and $Y$ are separated by $m$, then the same proof holds provided that the quantities $a_0, a_1, a_2, b_1, b_2, b_3$ are signed.  \hfill Q.E.D.
\newline

The \emph{(signed) moment} of the finite point-mass system $\mathcal{X} = \{(X_i, x_i), i = 1, 2, 3, ...,n\}$  with respect to the \underline{directed} straight line $m$ is
\[M_m(\mathcal{X}) = \sum_{i=1}^n\sigma_m(X_i)M_m(X_i,x_i)
\] where $\sigma_m(X)=  1, -1, 0$ according as $X$ is in the left half-plane of $m$, right half-plane of $m$ or on $m$ itself. The finite point-mass system $\mathcal{X}$ is said to be $balanced$ with respect to the directed straight line $m$ provided
\[M_m(\mathcal{X}) = 0.
\]It is clear that if $m$ and $m'$ are reverses of each other, then for every finite system  $\mathcal{X}$ we have 
\[M_m(\mathcal{X}) = - M_{m'}(\mathcal{X})
\]and \[M_m(\mathcal{X}) = 0 \ \ \ \ \emph{\emph{if and only if}}\ \ \ \ M_{m'}(\mathcal{X}) = 0
\]

\begin{corollary}For every finite point-mass system $\mathcal{X}$ and directed straight line $m$ 
\[M_m(\mathcal{X}) = M_m(\mathcal{C(X)})
\]
\end{corollary}PROOF: This follows from Theorem 3.4 by  induction.\hfill Q.E.D.

\begin{corollary}Every finite point-mass system is balanced with respect to every straight line that contains its centroid. \end{corollary}
Proof: This follows immediately from Corollary 3.5. \hfill Q.E.D.
\newline

\begin{figure}\center\includegraphics{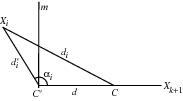}\caption{}\end{figure}
We now generalize Eq'n (3) to a formula for the mass of an arbitrary finite point-mass system. 
\begin{theorem}Let $\mathcal{X} = \{(X_i, x_i), i = 1, 2, 3, ...,n\}$ be a finite point-mass system and let
\[(C, c) = \mathcal{C}(\mathcal{X}).
\]If $d_i$ denotes the hyperbolic distance from $X_i$ to $C$ for $i = 1, 2, ..., n$, then
\[c = \sum_{i=1}^nx_i\cosh d_i.
\]
\end{theorem}PROOF: We proceed by mathematical induction on $n$.  The case $n = 1$ is self evident.  The case $n = 2$ is Eq'n (3).  The case $n=3$ is the last paragraph of the proof of Proposition 3.3.  Assume the theorem has been proved for $n = k$ and let
\[(C, c) = \mathcal{C}(\mathcal{X})\  \emph{\emph{where}}\  \mathcal{X} = \{(X_i, x_i), i = 1, 2, 3, ...,k+1\}\]  
and
\[(C', c') = \mathcal{C}(\mathcal{X'})\  \emph{\emph{where}}\  \mathcal{X'} = \{(X_i, x_i), i = 1, 2, 3, ...,k\}.\] 
Let (see Fig. 7) 
\[ d_i = d(X_i, C)\ \ i = 1, 2, ..., k+1  \ \ \ \ d'_i = d(X_i, C')\ \ i = 1, 2, ..., k\]
and\[d = d(C, C').\]
By the induction hypothesis
\[x' = \sum_{i=1}^kx_i\cosh d_i'  \]
Since
\[(C, c) = (C',c')*(X_{k+1},x_{k+1})
\]it follows from Eq'n (4) and the Law of Cosines that

\[x = x'\cosh d+ x_{k+1}\cosh d_{k+1} 
\]

\[= \sum_{i=1}^kx_i\cosh d_i'\cosh d+ x_{k+1}\cosh d_{k+1} \]
\[= \sum_{i=1}^kx_i(\cos\alpha_i\sinh d_i'\sinh d + \cosh d_i)+x_{k+1}\cosh d_{k+1}
\]
\[=-\left(\sum_{i=1}^kx_i(\sin(\alpha_i-\pi/2)\sinh d_i'\right)\sinh d + \sum_{i=1}^{k+1}x_i\cosh d_{i}
\]which, by Corollary 3.6, where the line $m$ in question passes through $C'$ and is perpendicular to $CC'$, equals
\[-0\sinh d +\sum_{i=1}^{k+1}x_i\cosh d_i\]
\[= \sum_{i=1}^{k+1}x_i\cosh d_i.
\] 
\hfill Q.E.D.
\newline

\begin{proposition}Two finite point-mass systems have the same moment with respect to a straight line if and only if it contains at least one of the centroids of their centroids.
\end{proposition}
PROOF:  This follows immediately from Proposition 3.2 and Corollary 3.5.\hfill Q.E.D.

\section{CENTROIDS OF LAMINAE}

A $region$ is a compact  subset of the hyperbolic plane of finite positive measure. A $lamina$ $\mathcal{L}$ is a pair $(L,\lambda)$ where $L$ is a region and $\lambda$ is a continuous non-negative valued function on $L$ such that
\[\int\int_L\lambda(X)dA > 0.
\] The value $\lambda(X)$ is the \emph{density} of $\mathcal{L}$ at $X$.  The lamina is said to be $uniform$ if its density is constant throughout $L$. The maximum value of $\lambda$ over $L$ is denoted by $\Lambda(\mathcal{L})$.  If  $P$ is any point and $p$  is any straight line, then
\[\Gamma_p(L) = \max_{X\in L}\{\cosh[d(X,p)]\}\ \ \ \ \emph{\emph{and}}\ \ \ \ \ \Gamma_P(L) = \max_{X\in L}\{\cosh[d(X,P)]\}
\]

 A $decomposition$ of $\mathcal{L}$ is  a family of sets $\tilde{L} = \{L_1, L_2, ..., L_n\} $ such that
\[L = L_1 \cup L_2 \cup \cdot\cdot\cdot \cup L_n\] 
where distinct $L_i$'s intersect in sets of measure 0. If each of the $L_i$'s has diameter less than $\delta$, this is a $\delta$-$decomposition$.  A $\delta$-$transversal$ of the $\delta$-decomposition $\tilde{L}$ is a point-mass system $\mathcal{X}$ = $\{(X_1, x_1), $ $(X_2 , x_2), ..., (X_n, x_n)\}$ such that
\[X_i \in L_i \ \  \emph{\emph{and}}\ \ x_i = \lambda(X_i)\emph{\emph{area}}(L_i), \ \ \ \ i = 1, 2, ..., n.
\]

Let  $m$ be a directed straight line.  We define the \emph{moment of} $\mathcal{L}$  \emph{with respect to m} as

\[M_m(\mathcal{L}) = \int\int_L\sigma_m(X)\lambda(X)\sinh[ d(X,m)]dA\]
where $dA$ is the area element.  The following technical lemma is needed for the proof of the crucial Theorem 4.2.

\begin{lemma}Let the distinct straight lines $m_1$ and $m_2$ intersect in the point $P$.  For each point $X\neq P$ and $i = 1, 2, $ let $\alpha_i(X)$ denote the non-obtuse angle between $m_i$ and $XP$, and let  $\beta$ be one of the angles determined by $m_1$ and $m_2$.  Then there exists a positive number $\Delta(m_1,m_2)$ such that
\[\min \left\{\csc \left(\alpha_1(X)\right),\csc \left(\alpha_2(X)\right)\right\} < \Delta(m_1, m_2)\ \ \ \ \emph{\emph{for all}}\  X \neq P.
\]\end{lemma}
PROOF:  Suppose, by way of contradiction, that  $\Delta(m_1, m_2)$ does not exist.  It follows that for each positive integer $n$ there exists a point $X_n \neq P$ such that
\[\csc \left(\alpha_1(X_n)\right) \geq n\ \ \ \emph{\emph{and}}\ \ \ \csc \left(\alpha_2(X_n)\right) \geq n.
\]It follows that
\begin{equation}\lim_{n\rightarrow\infty} \alpha_i(X_n) = 0,\ \ \ \  i = 1, 2.
\end{equation}However, it is clear from Figure 8a that for each $X\neq P$ either
\[\alpha_1(X) + \alpha_2(X) \in \{ \beta, \pi-\beta\}\]or\[|\alpha_1(X) - \alpha_2(X)| \in  \{ \beta, \pi-\beta\}.
\]Since the lines $m_1, m_2$ are distinct it follows that the angles $\beta$ and $\pi-\beta$ are neither 0 nor $\pi$ so that Eq'n (6) above leads to a contradiction.  Hence the required $\Delta(m_1, m_2)$ exists.\hfill Q.E.D.
\newline

\begin{theorem} Let $\mathcal{L} = (L, \lambda)$ be a lamina and let $m_i, i = 1, 2, 3$ be three concurrent straight lines such that 
\[M_{m_1}(\mathcal{L}) = M_{m_2}(\mathcal{L}) = 0.
\]Then
\[M_{m_3}(\mathcal{L}) = 0.\]
\end{theorem}
PROOF: Let $P$ be the intersection of all the $m_i$'s and suppose, by way of contradiction, that 
\[a = \left|M_{m_3}(\mathcal{L})\right| \neq 0.\] 
\begin{figure}\center\includegraphics{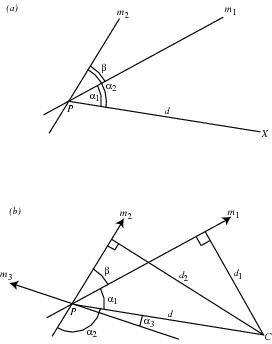}\caption{}\end{figure}
Let $n$ be an integer greater than $1 + \Delta(m_1, m_2)$.  By the definition of integrals, there exists a partition $\tilde{L}$ of $L$ and a transversal $\mathcal{X}$ of $\tilde{L}$ such that
\begin{equation}\left|M_{m_i}(\mathcal{X})\right| = \left|M_{m_i}(\mathcal{X}) - M_{m_i}(\mathcal{L})\right| < \frac{a}{n}    \ \ \ \ \  i = 1, 2                           
\end{equation}and
\begin{equation}\left|M_{m_3}(\mathcal{X}) - M_{m_3}(\mathcal{L})\right| < \frac{a}{n}.
\end{equation} Direct the $m_i$'s so that
\[M_{m_i}(\mathcal{X}) \geq 0,\ \ \ \  i = 1, 2, 3.
\]
Let $\mathcal{C(X)} = (C, c)$, $d_i = d(C,m_i)$ and let $\alpha_i$ be
 either of the positive angles between $m_i$ and $CP$ (see Fig. 8b).  Then,  for i = 1, 2, 3, 
\[M_{m_i}(\mathcal{X})= M_{m_i}(\mathcal{C(X)}) = c \sinh d_i = c \sinh d \sin \alpha_i
\]so that
\[\frac{M_{m_1}(\mathcal{X})}{ \sin \alpha_1} = \frac{M_{m_2}(\mathcal{X})}{ \sin \alpha_2 }= \frac{M_{m_3}(\mathcal{X})}{ \sin \alpha_3}  = c \sinh d.
\]It follows from Eq'n (8) that
\[-\frac{a}{n} < M_{m_3}(\mathcal{X}) - a 
\]or
\[a\left(1 - \frac{1}{n}\right) < M_{m_3}(\mathcal{X}).
\]
Eq'n (7) yields
\[
M_{m_i}(\mathcal{X}) < \frac{a}{n} ,\ \ \ i = 1, 2.
\]Hence
\[\csc\alpha_i \geq \frac{\sin\alpha_3}{\sin \alpha_i} = \frac{M_{m_3}(\mathcal{X})}{M_{m_i}(\mathcal{X})}  > n - 1 > \Delta(m_1,m_2),\ \ \ i = 1, 2,
\]which contradicts the definition of $\Delta(m_1,m_2)$. \hfill Q.E.D.\newline

It follows that for any  lamina $\mathcal{L} = (L,\lambda)$, the straight lines with respect to which $\mathcal{L}$ is balanced (i.e., has moment 0) are concurrent and this common point is the \emph{location of the} \emph{center of mass} of $\mathcal{L}$.  If this location is denoted by $C(\mathcal{L})$ then, consistently with Theorem 3.7, the \emph{mass} of $\mathcal{L}$ is defined as

\[ c(\mathcal{L}) = \int\int_L\lambda(X)\cosh[d(X, C(\mathcal{L})]dA.
\]
The pair $(C(\mathcal{L})), c(\mathcal{L}))$ is the \emph{center of mass} or \emph{centroid} $\mathcal{C(L)}$ of $\mathcal{L}$.

\begin{proposition}Let $\mathcal{L} = (L,\lambda)$ be a lamina.  Then for every $\epsilon > 0$ there is a $\delta > 0$ such that for every $\delta$-transversal $\mathcal{X}$ of $\mathcal{L}$
\begin{equation}d(C(\mathcal{X)},C(\mathcal{L})) < \epsilon\ \ \ \ \emph{and}\ \ \ \ |c(\mathcal{X})-c(\mathcal{L})| < \epsilon,
\end{equation}where $\mathcal{C(L)} = (C(\mathcal{L}), c(\mathcal{L}))$   and  $\mathcal{C(X)} = (C(\mathcal{X}), c(\mathcal{X}))$.
\end{proposition}
PROOF:  Let $\epsilon > 0$ and let $\delta$ be such that for every  $\delta$-decomposition $\tilde{L} = \{L_i, i = 1, 2, ..., n\}$ of $\mathcal{L} $ and for all directed straight lines $m$ through $C(\mathcal{L})$

\[0 < \max_{X\in L_i}\{\lambda(X)\} - \min_{X\in L_i}\{\lambda(X)\}< \frac{\epsilon}{2\Gamma_{C(\mathcal{L})}(L) \emph{\emph{area}}(L)}
\]
\[0 < \max_{X\in L_i}\{\sinh[d(X,m)]\} - \min_{X\in L_i}\{\sinh[d(X,m)]\}< \frac{\epsilon}{2\Lambda(L) \emph{\emph{area}}(L)}\]for all $i = 1, 2, ..., n$.

 Then, for every $m$ through $\mathcal{C(L)}$
\[|M_m(\mathcal{X})| = \left|M_m(\mathcal{X}) - M_m(\mathcal{L})\right|\]\[ = 
\left|\sum_{i=1}^n\sigma_m(X_i)\lambda(X_i)\sinh[d(X_i,m)]\emph{\emph{area}}(L_i) \right.\]\[ \left.-\int\int_L\sigma_m(X)\lambda(X)\sinh[d(X,m)]dA\right|
\]
\[\leq \sum_{i=1}^n \left[\max_{X\in L_i}\{\lambda(X)\sinh[ d(X,m)] \}\right.\]\[\left. -\min_{X\in L_i}\{\lambda(X)\sinh[ d(X,m)] \}\right]\emph{\emph{area}}(L_i)
\]
\[\leq  \sum_{i=1}^n\left[\max_{X\in L_i}\{\lambda (X)\} \max_{X\in L_i}\{\sinh[d(X_i,m)]\}\right.
\]
\[- \left.\min_{X\in L_i}\{\lambda (X)\} \min_{X\in L_i}\{\sinh[d(X_i,m)]\}\right]\emph{\emph{area}}(L_i)
\]
\[=  \sum_{i=1}^n\left[\max_{X\in L_i}\{\lambda (X)\}\left( \max_{X\in L_i}\{\sinh[d(X_i,m)]\}-\min_{X\in L_i}\{\sinh[d(X_i,m)]\}\right)\right.
\]
\[+ \left(\left.\max_{X\in L_i}\{\lambda (X)\}-\min_{X\in L_i}\{\lambda (X)\}\right) \min_{X\in L_i}\{\sinh[d(X_i,m)]\}\right]\emph{\emph{area}}(L_i)
\]
\[< \sum_{i=1}^n\left[ \Lambda(L) \frac{\epsilon}{2\Lambda(L)\cdot 
\emph{\emph{area}}(L)}\emph{\emph{area}}(L_i)\right.
\]\[\left.+ \frac{\epsilon}{2\Gamma_{C(\mathcal{L})}(L)\emph{\emph{area}}(L)}\Gamma_{C(\mathcal{L})}(L) \emph{\emph{area}}(L_i)
\right]
\]
\[<\frac{\epsilon}{2} + \frac{\epsilon}{2} = \epsilon.\]

Suppose, by way of contradiction, that the first inequality of (9) is false.  Then there exists an $\epsilon_0 > 0$ such that for every positive integer $k$ there is a $\frac{1}{k}$-transversal $\mathcal{X}^{(k)}$ such that 
\[\sum_{i=1}^{n_k}\lambda(X_i)\emph{\emph{area}}(L_i) \geq \frac{1}{2}\int\int_L\lambda(X)dA
\]
and
\[ \sinh[d(C(\mathcal{X}^{(k)}),C(\mathcal{L}))] \geq d(C(\mathcal{X}^{(k)}),C(\mathcal{L})) \geq \epsilon_0.\]\newline
However, by the first part of the proof and Corollary 3.5, for all sufficiently large $k$ and for that $m$ that is perpendicular to the straight line joining $C(\mathcal{X}^{(k)})$ to $C(\mathcal{L})$

\[c(\mathcal{X}^{(k)})\epsilon_0 \leq c(\mathcal{X}^{(k)})\sinh[d(C(\mathcal{X}^{(k)}),C(\mathcal{L})]  \]
\[=c(\mathcal{X}^{(k)})\sinh[d(C(\mathcal{X}^{(k)}),m)]  \]
\[= |M_m(\mathcal{X}^{(k)})| 
< \frac{\epsilon_0}{2}\int\int_L\lambda(X)dA.
\]Hence
 \[\frac{1}{2}\int\int_L\lambda(X)dA > c(\mathcal{X}^{(k)}) \]
 \[= \sum_{i=1}^{n_k}\lambda(X_i)\emph{\emph{area}}(L_i)\cosh[d(X_i,C(\mathcal{X}^{(k)})) \]
 \[\geq \sum_{i=1}^{n_k}\lambda(X_i)\emph{\emph{area}}(L_i)  \geq \frac{1}{2}\int\int_L\lambda(X)dA\] 
 which is impossible.  
This establishes the first inequality of (9). 

 The second inequality now follows by standard  arguments. For any $\delta$-transversal $\mathcal{X}$ of $\mathcal{L}$ we have
\[|c(\mathcal{X}) - c(\mathcal{L})|\]\[ = \left|\sum_{i=1}^n\lambda(X_i)\emph{\emph{area}}(L_i)\cosh[d(X_i,C(\mathcal{X}))]
-\int\int_L\lambda(X)\cosh[d(X,C(\mathcal{L}))]dA\right|
\]
\[ \leq \left|\sum_{i=1}^n\lambda(X_i)\emph{\emph{area}}(L_i)\cosh[d(X_i,C(\mathcal{X}))]
- \sum_{i=1}^n\lambda(X_i)\emph{\emph{area}}(L_i)\cosh[d(X_i,C(\mathcal{L}))]\right|\]
\[+\left|\sum_{i=1}^n\lambda(X_i)\emph{\emph{area}}(L_i)\cosh[d(X_i,C(\mathcal{L}))] - \int\int_L\lambda(X)\cosh[d(X,C(\mathcal{L}))]dA\right|
\]

\[\leq \sum_{i=1}^n\lambda(X_i)\emph{\emph{area}}(L_i)\left|\cosh[d(X_i,C(\mathcal{X}))] - \cosh[d(X_i,C(\mathcal{L}))] \right|
\]
\[+\left|\sum_{i=1}^n\lambda(X_i)\emph{\emph{area}}(L_i)\cosh[d(X_i,C(\mathcal{L}))] - \int\int_L\lambda(X)\cosh[d(X,C(\mathcal{L}))]dA\right|
\]

\[\leq \sum_{i=1}^n\lambda(X_i)\emph{\emph{area}}(L_i)2\sinh[\emph{\emph{diameter}}(L)]\sinh[d(C(\mathcal{X}),C(\mathcal{L}))]
\]
\[+\left|\sum_{i=1}^n\lambda(X_i)\emph{\emph{area}}(L_i)\cosh[d(X_i,C(\mathcal{L}))] - \int\int_L\lambda(X)\cosh[d(X,C(\mathcal{L}))]dA\right|
\]
However, it is clear that each of the two summands of the above expression can be made arbitrarily small by choosing $\delta$ small enough.

\hfill Q.E.D.\newline

\begin{proposition}If $\mathcal{L}$ is a lamina, and $m$ is any straight line, then
\[M_m(\mathcal{L}) = M_m(\mathcal{C(L)}).
\]
\end{proposition}
PROOF: Let $m$ be a fixed straight line, let $\epsilon > 0$ be given and let $\mathcal{X}$ be a $\delta$-transversal of $\mathcal{L}$ such that
\[d(C(\mathcal{X}),C(\mathcal{L})) < \min\left\{\frac{\epsilon}{8c(\mathcal{L})\Gamma_m(L)},d(C(\mathcal{L}),m)\right\},\]

\[|c(\mathcal{X})-c(\mathcal{L})| < \min\left\{\frac{\epsilon}{4\Gamma_m(L)}, c(\mathcal{L})\right\}\]and  \[\left|M_m(\mathcal{X}) - M_m(\mathcal{L})\right| < \frac{\epsilon}{2}.
\]Set
\[\mathcal{C(X)} = (C_1, c_1),\ \ \ \mathcal{C(L)} = (C_2, c_2)\ \ \ \ d_i = d(C_i,m),\ \ \ \  i = 1, 2.
\]
Then $C(\mathcal{X})$ and $C(\mathcal{L})$ are on the same side of $m$ and \[\left|M_m(\mathcal{C(X)}) - M_m(\mathcal{C(L)}) \right| = |c_1\sinh d_1 - c_2\sinh d_2|
\]
\[\leq c_1|\sinh d_1 - \sinh d_2| + |c_1 - c_2|\sinh d_2
\]
\[\leq 2c_2 \Gamma_m(L)|d_1 - d_2| + \frac{\epsilon}{4\Gamma_m(L)}\Gamma_m(L)\]\[  \leq 2c_2 \Gamma_m(L) d(C_1, C_2) +\frac{\epsilon}{4}
\]\[<\frac{\epsilon}{4}+ \frac{\epsilon}{4} = \frac{\epsilon}{2}.
\]
It follows that
\[\left|M_m(\mathcal{L}) - M_m(\mathcal{C(L)})\right| \]
\[\leq \left|M_m(\mathcal{L}) - M_m(\mathcal{X)})\right|
+ \left|M_m(\mathcal{X}) - M_m(\mathcal{C(X)})\right|\]\[ + \left|M_m(\mathcal{C(X)}) - M_m(\mathcal{C(L)})\right|
\]
\[\leq \frac{\epsilon}{2} + 0 + \frac{\epsilon}{2} = \epsilon.
\]Since $\epsilon$ is arbitrary, the proposition follows. \hfill Q.E.D.

\begin{corollary}Two laminae have the same moment with respect to every directed straight line if and only if they have identical centroids.
\end{corollary}
PROOF:  This follows from Propositions 4.4 and 3.2. \hfill Q.E..D.

\begin{proposition}Let $\mathcal{L} = (L, \lambda)$ be a lamina, $\tilde{L} = \{L_1, L_2, ..., L_n\}$ a decomposition of $\mathcal{L}$ and set
\[\mathcal{L}_i = (L_i,\lambda_{|L_i}),\ \ \ \ \ \ i = 1, 2, ..., n.
\]Then
\[\mathcal{C(L)} = \mathcal{C(L}_1)*\mathcal{C(L}_2)* \cdot\cdot\cdot*\mathcal{C(L}_n).
\]
\end{proposition}
PROOF:  It follows from Proposition 4.4 and the additivity of integrals that for any directed straight line $m$
\[M_m[\mathcal{C(L)}] = M_m(\mathcal{L}) = M_m(\mathcal{L}_1) + M_m(\mathcal{L}_2) + \cdot\cdot \cdot + M_m(\mathcal{L}_n)
\]
\[= M_m(\mathcal{C(L}_1)) + M_m(\mathcal{C(L}_2)) + \cdot\cdot \cdot + M_m(\mathcal{C(L}_n))
\]
\[=M_m[\mathcal{C(L}_1)*\mathcal{C(L}_2)*\cdot\cdot\cdot*\mathcal{C(L}_3)].
\]
The validity of the proposition now follows from the arbitrariness of $m$ and Corollary 4.5. \hfill Q.E.D. \newline
\section{CENTROIDS OF LINEAR SETS}
We now briefly discuss the 1-dimensional analogs of laminae.
A \emph{linear set} $\mathcal{L} = (l, \lambda)$ is a non-empty, compact, and measurable subset  $l$ of a straight line in the hyperbolic plane, and a non-negative function $\lambda: l \rightarrow \mathcal{R}$ such that
\[\int_l\lambda(X)dX > 0.
\]  If $m$ is either of the directed straight lines that contain $l$ and $A$ is any point of $l$ then the \emph{moment of $\mathcal{L}$ with respect to $A$} is
\[M_A(\mathcal{L}) = \int_l\sigma_A(X)\lambda(X)\sinh[d(X,A)]dX.
\]where $\sigma_A(X)$ = 1 or -1 according as the direction from $A$ to $X$  agrees or disagrees with that of $m$, and $\sigma_A(A) = 0$.
The unique point $C$ of $m$ such that
\[M_C(\mathcal{L}) = 0
\]is the \emph{location of the centroid of } $ \mathcal{L}$.  In analogy with Theorem 3.7 and the definition of the mass of a lamina, the \emph{mass} of the linear set $\mathcal{L}$ is
\[\emph{\emph{c}}(\mathcal{L}) = \emph{\emph{mass}}(\mathcal{L}) = \int_l\lambda(X)\cosh[d(X,C)]dX
\]
 If $m$ is another directed straight line then the \emph{moment of $\mathcal{L}$ with respect to $m$} is
\[M_m(\mathcal{L}) = \int_l\sigma_m(X)\lambda(X)\sinh[d(X,m)]dX.
\]The pair
\[\mathcal{C(L)} = (C(\mathcal{L}),c(\mathcal{L}))
\]is the centroid of $\mathcal{L}$.
\begin{example}\end{example}The centroid of a hyperbolic line segment of length $d$ and uniform density 1 is located at its midpoint and its mass is defined to be
\[2\int_0^{d/2}\cosh x dx = 2\sinh (d/2).
\]

The following four propositions are linear analogs of Propositions 3.2 and 4.3 - 4.5 and their proofs, being simplifications of the 2-dimensional proofs are omitted.

\begin{proposition} If $\mathcal{L}$ is a linear set and $m$ is any straight line that contains $C(\mathcal{L})$, then $M_m(\mathcal{L}$) = 0.
\end{proposition}\hfill$\Box$\newline

\begin{proposition}Let $\mathcal{L} = (L,\lambda)$ be a linear set.  Then for every $\epsilon > 0$ there is a $\delta > 0$ such that for every $\delta$-transversal $\mathcal{X}$ of $\mathcal{L}$
\[d(C(\mathcal{X)},C(\mathcal{L})) < \epsilon\ \ \ \ \emph{\emph{and}}\ \ \ \ |c(\mathcal{X})-c(\mathcal{L})| < \epsilon,
\]where $\mathcal{C(L)} = (C(\mathcal{L}), c(\mathcal{L}))$   and  $\mathcal{C(X)} = (C(\mathcal{X}), c(\mathcal{X}))$.
\end{proposition}\hfill$\Box$\newline

\begin{proposition}If $\mathcal{L}$ is a linear set, and $m$ is any straight line, then
\[M_m(\mathcal{L}) = M_m(\mathcal{C(L)}).
\]
\end{proposition}\hfill$\Box$\newline

\begin{proposition}Let $\mathcal{L} = (L, \lambda)$ be a linear set, $\tilde{L} = \{L_1, L_2, ..., L_n\}$ a decomposition of $\mathcal{L}$ and set
\[\mathcal{L}_i = (L_i,\lambda_{|L_i}),\ \ \ \ \ \ i = 1, 2, ..., n.
\]Then
\[\mathcal{C(L)} = \mathcal{C(L}_1)*\mathcal{C(L}_2)* \cdot\cdot\cdot*\mathcal{C(L}_n).
\]
\end{proposition}\hfill$\Box$\newline

The following proposition is a mathematical  analog of Archimedes's "mechanical" method for finding volumes and centroids [Archimedes].

\begin{proposition}Let $\mathcal{L} = (L,\lambda)$ be a lamina, $\Pi$ a pencil of asymptotically parallel straight lines, and $m$ a straight line.  Suppose that for every $p \in \Pi$, the pair   $(L\cap p,\lambda_{|L\cap p})$ is a linear set whenever $L\cap p$ has positive 1-dimensional measure, and
\[M_m(L\cap p,\lambda_{|L\cap p}) = 0.\]  Then
\[M_m(\mathcal{L}) = 0.
\]
\end{proposition}
PROOF: We work in the upper half-plane model where
\[ds = \frac{\sqrt{x^2 + y^2}}{y}\ \ \ \ \emph{\emph{and}}\ \ \ \ dA = \frac{dxdy}{y^2}.
\]By symmetry it may be assumed that $\Pi$ consists of all the geodesics of the form
\[p_a = \{(a,y)\  |\  a\  \emph{\emph{is fixed and}}\  y > 0\}.
\]Let
\[N = \{x\ |\ L\cap p_x \neq \emptyset\}.
\]Then
\[M_m(\mathcal{L}) = \int\int_L\sigma_m(x,y)\lambda(x,y)\sinh d[((x,y),m)]\frac{dxdy}{y^2}
\]
\[= \int_N\left[\int_{L\cap p_x}\sigma_m(x,y)\lambda(x,y)\sinh d[((x,y),m)]\frac{dy}{y}\right]\frac{dx}{y}
\]
\[= \int_NM_m(L\cap p_x,\lambda_{|L\cap p_x)}\frac{dx}{y}
\]
\[= \int_N0\frac{dx}{dy} = 0.
\]\hfill Q.E.D.

\section{EXAMPLES}

The Euclidean analog of the following proposition [Ungar] is well known.
\begin{figure}\center\includegraphics{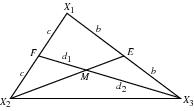}\caption{}\end{figure}
\begin{proposition}
Let $(C, c)$ =  $\mathcal{C}\{(X_1, w), (X_2, w), (X_3, w)\}$. Then $C$ is the point of intersection of the medians of $\Delta X_1X_2X_3$.\end{proposition}
Proof:  Let $E, F$ be the respective midpoints of the sides  $X_1X_3$ and $X_1X_2$ of  $\Delta X_1X_2X_3$ (Fig. 9). Then the centroid of $\{(X_1, w), (X_2, w)\}$ is the point-mass
\[(F, 2w\cosh c)
\]and hence the centroid of the system $\{(X_1, w), (X_2, w), (X_3, w)\}$ lies on the point $M$ of $X_3F$ such that
\[2w\cosh c \sinh d_1 = w \sinh d_2.
\]It follows that
\[\frac{\sinh X_1E}{\sinh E X_3}\frac{\sinh X_3M}{\sinh MF}\frac{\sinh F X_2}{\sinh X_2X_1} = \frac{\sinh b}{\sinh b}\frac{\sinh d_2}{\sinh d_1}\frac{\sinh c}{\sinh 2c}
\]\[= \frac{1}{1}\frac{2w \cosh c}{w}\frac{1}{2\cosh c} = 1.
\]Hence, by the converse to the theorem of Menelaus, the points $X_2, M,$ and $E$ are collinear.  Since the medians of the hyperbolic triangle are concurrent, their common intersection is also the location of the centre of mass in question.

\hfill Q.E.D.\newline

Some of the subsequent examples are worked out in a specific model that is based on a  general  geodesic polar  parametrization used by Gauss in [Gauss].  This Gaussian model presents the hyperbolic plane as a Riemannian geometry whose domain is the entire plane with polar coordinates $(\rho, \theta)$ and  metric [Gauss, Stahl]
\[d\rho^2 + \sinh^2\rho d\theta^2
\]  The geodesics of this metric are the Euclidean straight lines $\theta = c$ and the curves
\[\rho =\coth^{-1}(C\cos(\theta - \alpha)
\]where $\alpha$ is arbitrary and $C > 1$.  The area element of this metric is
\[dA = \sinh \rho d\rho d\theta.
\]
It is clear that mass is invariant under rigid motions and consequently the axes of reflections of a region contain its centroid.  In particular the centroid of a uniform disk is located at its center.

\begin{proposition}The mass of a disk  of uniform density 1 and hyperbolic radius r is
\[\pi\sinh^2r.
\]\end{proposition}
PROOF: We employ the Gauss model and assume that the disk is centered at the origin which coincides with its centroid.  By the definition of mass, the mass of this disk is 
\[\int_0^{2\pi}\int_0^r\cosh \rho \sinh \rho d\rho d\theta = \pi \sinh^2 r.
\]\hfill Q.E.D.\newline

This formula is particularly interesting for the following reason.  As was noted above, many hyperbolic formulas can be obtained from their Euclidean analogs by the heuristic means of replacing a certain length  $d$  by  $\sinh d$.  One of the exceptions to this informal rule is the area of a circle of radius $r$.  The Euclidean formula is 
\[\pi r^2
\]whereas the hypebolic formula is
\[4\pi\sinh^2\left(\frac{r}{2}\right).
\]Thus, it would seem that while in Euclidean geometry area and uniform mass are essentially equivalent, in hyperbolic geometry, where they are distinct, sometimes it is the notion of mass that is better behaved (by Euclidean standards, of course).  Another instance is offered in Proposition 6.5.

We next turn to some uniform wedges; first their centroids are located and then their masses are computed. Let  $D_n = D_n(r)$ denote the lamina consisting of the  subset 

\[\{(\rho, \theta) \in D_n(r) \ | \  -\frac{\pi}{n} \leq \theta \leq \frac{\pi}{n}\}
\]of the disk $D(r)$ with uniform density 1 (Fig. 10).  Let  $d_n$ denote the distance from the origin   $O$  to $C(D_n)$ and  $R = R_{O,2\pi/n}$ denote the counterclockwise rotation by the angle $2\pi/n$ about  $O$.    Then, by symmetry, Proposition 4.6,  and the Law of Cosines

\begin{equation}\pi\sinh^2r = \emph{\emph{mass}}(D(r)) = \sum_{i=1}^n\emph{\emph{mass}}(R^i(D))d_n = n\ \emph{\emph{mass}}(D_n)\cosh d_n
\end{equation}
\[= n\cosh d_n\int_{-\pi/n}^{\pi/n}\int_0^r\cosh[d (C(D_n),X)] \sinh \rho d\rho d\theta
\] 
\[= n\cosh d_n\int_{-\pi/n}^{\pi/n}\int_0^r(\cosh d_n \cosh \rho - \cos\theta \sinh d_n \sinh \rho)\sinh \rho d\rho d\theta
\]
\[= n \cosh^2 d_n\int_{-\pi/n}^{\pi/n}\int_0^r\cosh \rho \sinh \rho d\rho d\theta\]\[ - n\  \cosh d_n\sinh d_n\int_{-\pi/n}^{\pi/n}\int_0^r\cos\theta \sinh^2 \rho d\rho d\theta
\]
\[=\cosh^2 d_n\pi\sinh^2 r - n\cosh d_n\sinh d_n\cdot2\sin\frac{\pi}{n}\int_0^r\frac{\cosh 2\rho - 1}{2}d\rho
\]
\[=\pi\cosh^2 d_n\sinh^2 r - n\sinh d_n\cosh d_n\sin \frac{\pi}{n}\left(\frac{\sinh 2r}{2} - r\right).
\]

Dvision by $\pi\sinh^2r$ yields \begin{figure}
\center
\includegraphics{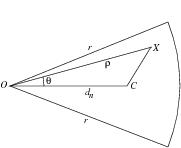}
\caption{}
\end{figure}

\[1 = \cosh^2 d_n - \frac{n}{\pi}\sin \frac{\pi}{n}\sinh d_n\cosh d_n\  \frac{\sinh 2r - 2r}{\cosh 2r - 1}\]
or \[\sinh^2d_n = \frac{n}{\pi}\sin\frac{\pi}{n}\sinh d_n\cosh d_n\ \frac{\sinh 2r - 2r}{\cosh 2r - 1}
\]or
\[\tanh d_n = \frac{n}{\pi}\sin \left(\frac{\pi}{n}\right) \frac{\sinh 2r - 2r}{\cosh 2r - 1}.
\]
It follows that
\[\frac{d_n}{r} = \frac{n}{\pi}\sin\left(\frac{\pi}{n}\right)\ \left( \frac{2}{3} + O(r^2)\right)
\]in comparison to the Euclidean analog of 2/3.
\[\frac{2n}{3\pi}\sin\left(\frac{\pi}{n}\right)
\]

It follows from Eq'n (10) that
\[\emph{\emph{mass}}(D_n) = \frac{\pi\sinh^2r}{n\cosh d_n}
\]
from which is obtained
\[\emph{\emph{mass}}\left(D_n(r)\right) = \frac{\pi\sinh^2r}{n}\sqrt{1 - \left[\frac{n}{\pi}\sin\left(\frac{\pi}{n}\right)\frac{\sinh 2r - 2r}{\cosh 2r -1}\right]^2}.
\]

We next turn to the centroid of the uniform triangular lamina.  A technical lemma sets the stage for a short proof that makes use of a mathematical analog of the "mechanical method" of Archimedes.  The Euclidean centroid of the uniform triangle lamina is, of course, well known. 

\begin{lemma}Let $\Delta ABC$ be a hyperbolic triangle with  points D, E, F, on the respective sides AB, BC, AC, such  that EF is asymptotically parallel to BC and let  $G = AD \cap EF$.  Then
\[\frac{\sinh BD}{\sinh CD} = \frac{\sinh GF}{\sinh GE}
\]

\end{lemma}
PROOF: In Figure 11, apply the Theorem of Menelaus to $\Delta FBH$ twice to obtain
\[\frac{\sinh  FA}{\sinh  AB}\ \frac{\sinh  BD}{\sinh  DH} \ \frac{\sinh  HG'}{\sinh  G'F} = \frac{\sinh  FA}{\sinh  AB}\ \frac{\sinh  BC}{\sinh  CH} \ \frac{\sinh  HE'}{\sinh  E'F}
\]Two similar applications to $\Delta E'CH$ yield
\[\frac{\sinh  E'A}{\sinh  AC}\ \frac{\sinh  CB}{\sinh  BH} \ \frac{\sinh  HF}{\sinh  FE'} = \frac{\sinh  E'A}{\sinh  AC}\ \frac{\sinh  CD}{\sinh  DH} \ \frac{\sinh  HG'}{\sinh  G'E'}
\]The multiplication of these two equations simplifies to
\[\frac{\sinh BD}{\sinh G'F}\ \frac{\sinh HF}{\sinh BH} = \frac{\sinh HE'}{\sinh CH}\ \frac{\sinh CD}{\sinh G'E'}
\]or
\[\frac{\sinh BD}{\sinh CD}\ \frac{\sinh HF}{\sinh HE'} = \frac{\sinh G'F}{\sinh G'E'}\ \frac{\sinh BH}{\sinh CH}
\]Since the limiting position of $E'$ and $G'$ as $H$ recedes to infinity along $BC$ (and $F$ is held fixed)  are $E$ and  $G$, respectively, it follows that 
\[\frac{\sinh BD}{\sinh CD} = \frac{\sinh GF}{\sinh GE'}
\]\hfill Q.E.D.
\begin{figure}
\center
\includegraphics{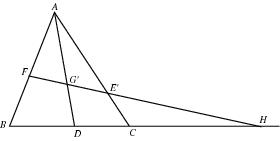}
\caption{}
\end{figure}

\begin{theorem}
The center of mass of a uniform triangle is located at the intersection of its medians.
\end{theorem}
PROOF:  It suffices to show that the uniform triangle is balanced with respect to its medians.  Let $\Delta ABC$ be such a triangle and $AD$ its median (Fig. 12).  By Lemma 5.4, if $FE$ is asymptotically parallel to $BC$, then
\[FG = EG.\]Hence,
\[M_{AD}(FG) = M_{AD}(EG)\]
and so, by Proposition 5.6.
\[M_{AD}(\Delta ABD) = M_{AD}(\Delta ACD).
\]\hfill Q.E.D.\newline 
\begin{figure}
\center
\includegraphics{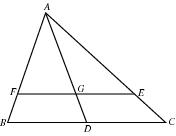}
\caption{}
\end{figure}

In both the statement and the proof below, the index $i$ is computed modulo 3.

\begin{proposition} Let $\Delta X_1X_2X_3$ be a triangular lamina with uniform density 1.  Then

\[\emph{mass}(\Delta X_1X_2X_3) = \frac{1}{2}\sum_{i=1}^3 \sinh [(d(O,X_iX_{i+1})]d(X_i,X_{i+1}).\]
\end{proposition}
PROOF:  To find the mass of the triangle we may assume that $O$, the intersection point of the medians, is the origin of a Gaussian parametrization of the hyperbolic plane (Fig. 13).  Then
	
\[\emph{\emph{mass}}(\Delta X_1X_2X_3) = \int\int_{\Delta X_1X_2X_3}\cosh \rho dA\]
\[ = \sum_{i=1}^3 \int\int_{\Delta OX_iX_{i+1}}\cosh \rho dA.
\]

Let
\[\rho_i = \rho_i(\theta) = \coth^{-1}\left(C_i\cos(\theta - \alpha_i)\right)
\]be the equation of the geodesic joining $X_{i+1}$ and $X_{i+2}$.  If, for $i = 1, 2, 3,$ $\theta_i$ is the angle from the horizontal axis to the geodesic $OX_i$  then
\[\int\int_{\Delta OX_iX_{i+1}}\cosh \rho dA = \int_{\theta_i}^{\theta_{i+1}}\int_0^{\rho_{i+2}(\theta)}\cosh\rho\sinh\rho d\rho d\theta
\]

\[=\frac{1}{2}\int_{\theta_i}^{\theta_{i+1}}\sinh^2\rho_{i+2}d\theta
\]
\[=\frac{1}{2}\int_{\theta_i}^{\theta_{i+1}}[\sinh (\coth^{-1}(C_{i+2} \cos (\theta - \alpha_{i+2})]^2d\theta
\]
\[=\frac{1}{2}\int_{\theta_i}^{\theta_{i+1}}\frac{d\theta}{C_{i+2}^2\cos^2(\theta - \alpha_{i+2})-1}.
\]

\begin{figure}\center\includegraphics{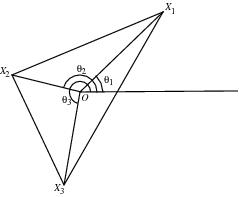}\caption{}\end{figure}

On the other hand, the length of the geodesic segment joining $X_iX_{i+1}$ is
\[d(X_i,X_{i+1}) = 
\int_{\theta_i}^{\theta_{i+1}}\sqrt{d\rho_{i+2}^2 + \sinh^2\rho_{i+2} d\theta^2}\]
\[ = \int_{\theta_i}^{\theta_{i+1}}\sqrt{\frac{C_{i+2}^2\sin^2(\theta - \alpha_{i+2})}{(C_i^2\cos^2(\theta - \alpha_{i+2}) - 1)^2} +  \frac{1}{C_{i+2}^2\cos^2(\theta - \alpha_{i+2}) - 1} }d\theta
\]
\[= \sqrt{C_{i+2}^2 - 1}\int_{\theta_i}^{\theta_{i+1}}\frac{d\theta}{C_{i+2}^2\cos^2(\theta - \alpha_{i+2})-1} \]
\[= \sqrt{C_{i+2}^2 - 1}\int\int_{\Delta OX_iX_{i+1}}\cosh \rho dA
\]

Set $d_i = d(O,X_{i+1}X_{i+2})$.  Then
\[\sqrt{C_i^2 - 1} = \sqrt{\coth^2d_i - 1} = \emph{\emph{csch}}\ d_i
\]
Hence,

\[\int\int_{\Delta OX_iX_{i+1}}\cosh \rho dA = \frac{\sinh d_{i+2}}{2}d(X_i,X_{i+1})
\]
and the  proposition now follows immediately. \hfill Q.E.D.
\newline

The Euclidean analog of our last proposition is also well known.
\begin{proposition}The mass of the regular n-gon of in-radius $r$ is half the product of its perimeter with $\sinh r$.
\end{proposition}
\begin{figure}\center\includegraphics{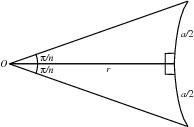}\caption{}\end{figure}
PROOF: Once again we work in the Gauss model of the hyperbolic plane.  Set $C = \coth r$ and let $a$ be the hyperbolic length of one of the polygon's sides (see Fig. 14).  Then one side of the polygon  is parametrized as 
\[ \rho = \coth^{-1}(C \cos\theta),\ \ \ \ -\pi/n \leq \theta \leq \pi/n.
\]  It follows from  the symmetry of the polygon that its mass equals
\[2n\int_0^{\pi/n}\int_0^{\coth^{-1}(C\cos\theta)}\cosh \rho \sinh \rho d\rho d\theta
\]
\[= n\int_0^{\pi/n}[\sinh (\coth^{-1}(C \cos \theta)]^2d\theta
\]
\[= n \int_0^{\pi/n}\left[\frac{\sqrt{\frac{C\cos\theta + 1}{C\cos\theta - 1}} - 
\sqrt{\frac{C\cos\theta - 1}{C\cos\theta + 1}}}{2}\right]^2d\theta
\]
\[= n \int_0^{\pi/n}\frac{d\theta}{C^2\cos^2\theta - 1} = \frac{n}{\sqrt{C^2-1}}\ \tanh^{-1}\left[\frac{\tan (\pi/n)}{\sqrt{C^2-1}}\right]\]
\[=n\ \sinh r \tanh^{-1}\left[\tan(\pi/n)\sinh r\right] \]\[= n\ \sinh r \tanh^{-1}\left[\tanh \left(\frac{a}{2}\right)\right] = \frac{na\sinh r}{2}.
\]
\hfill Q.E.D.

The area of the above regular polygon is well known to be
\[(n-2)\pi - 2n\beta.  
\] Thus the mass of the uniform regular polygon is also "better behaved" than its area.

\section{Appendix}

\begin{figure}\center\includegraphics{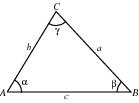}\caption{}\end{figure}

\begin{theorem}Let $\Delta ABC$ be the hyperbolic triangle of Figure 15.  Then
\[\frac{\sinh a}{\sin \alpha} = \frac{\sinh b}{\sin \beta} = \frac{\sinh c}{\sin \gamma}\ \ \ \ \ \ \emph{\emph{(Law of Sines)}}
\]
\[\cosh a = \cosh b \cosh c - \cos\alpha\sinh b \sinh c \ \ \ \ \ \ \emph{\emph{(Law of Cosines)}}
\]
\end{theorem}
\begin{theorem} Let P, Q, R,  be points on the respective extended sides  AB, BC, AC  of  the hyperbolic $\Delta ABC$.  Then
\newline\newline
Theorem of Ceva:
\[AP, BQ, CR\   \emph{are concurrent}\] 

if and only if \[\frac{\sinh AR}{\sinh RB} \frac{\sinh BP}{\sinh PC}\frac{\sinh CQ}{\sinh QA} = 1;
\]Theorem of Menelaus:
\[P, Q, R\   \emph{are collinear}\] 

if and only if \[\frac{\sinh AR}{\sinh RB} \frac{\sinh BP}{\sinh PC}\frac{\sinh CQ}{\sinh QA} = -1.
\]
\end{theorem}
\Large\textbf{Acknowledgements:} \normalsize The author is indebted to his colleague David Lerner for his help and patience.

\end{document}